\documentclass[12pt]{amsart}

\usepackage[all]{xy}
\usepackage{amscd}
\usepackage{amssymb}
\usepackage{calrsfs}

\newtheorem{thm}{Theorem}[section]
\newtheorem{lem}[thm]{Lemma}

\newtheorem{prop}[thm]{Proposition}
\newtheorem{conj}[thm]{Conjecture}

\theoremstyle{definition}

\newtheorem{rem}[thm]{Remark}
\newtheorem{claim}[thm]{Claim}
\newtheorem{defn}[thm]{Definition}

\newtheorem{ex}[thm]{Example}

\theoremstyle{remark}

\numberwithin{equation}{section}

\def\F{{\mathbb F}}

\def\Q{{\mathbb Q}}
\def\Z{{\mathbb Z}}
\def\C{{\mathbb C}}
\def\O{{\mathcal O}}
\def\P{{\mathbb P}}

\def\X{{\mathfrak X}}

\def\Fr{\text{\rm Fr}}
\def\Gr{\text{\rm Gr}}

\def\Gal{\text{\rm Gal}}

\def\Im{\text{\rm Im}\,}
\def\Ker{\text{\rm Ker}\,}

\def\Pic{\text{\rm Pic}}
\def\PGL{\text{\rm PGL}}

\def\bb{\hspace*{-2.5mm}}

\begin{document}

\title[Weight-monodromy conjecture for threefolds]
{Weight-monodromy conjecture for certain threefolds
in mixed characteristic}

\author[Tetsushi Ito]{Tetsushi Ito}
\address{Department of Mathematical Sciences,
University of Tokyo, 3-8-1 Komaba, Meguro, Tokyo 153-8914, Japan}
\email{itote2\char`\@ms.u-tokyo.ac.jp}

\address{Max-Planck-Institut F\"ur Mathematik,
Vivatsgasse 7, D-53111 Bonn, Germany}
\email{tetsushi\char`\@mpim-bonn.mpg.de}

\date{\today}

\begin{abstract}
The weight-monodromy conjecture claims the coincidence
of the weight filtration and the monodromy filtration,
up to shift, on the $l$-adic \'etale cohomology of
a proper smooth variety over a complete discrete valuation field.
Although it has been proved in some cases,
the case of dimension $\geq 3$ in mixed characteristic
is still open so far.
The aim of this paper is to give a proof of
the weight-monodromy conjecture
for a threefold which has a projective strictly semistable model
such that, for each irreducible component of the special fiber,
the Picard number is equal to the second $l$-adic Betti number.
Our proof is based on a careful analysis of
the weight spectral sequence of Rapoport-Zink
by the Hodge index theorem for surfaces.
We also prove a $p$-adic analogue by using
the weight spectral sequence of Mokrane.
\end{abstract}

\maketitle

\section{Introduction}
\label{SectionIntroduction}

Let $K$ be a complete discrete valuation field with
residue field $F$,
$\O_K$ the ring of integers of $K$,
and $l$ a prime number different from the characteristic of $F$.
Let $X$ be a proper smooth variety over $K$ and
$V = H^w(X_{\overline{K}},\Q_l)$ the $l$-adic
\'etale cohomology group
of $X_{\overline{K}} = X \otimes_{K} \overline{K}$ on which
the absolute Galois group $\Gal(\overline{K}/K)$ acts.
Let $M$ be the monodromy filtration on $V$,
and $W$ the weight filtration on $V$.
$M$ is defined by the action of the inertia group $I_K$.
If $F$ is finite, $W$ is defined by Frobenius eigenvalues.
If $X$ has a proper strictly semistable model over $\O_K$,
$W$ can also be defined by the weight spectral sequence
of Rapoport-Zink (for details, see
\S \ref{SectionFiltration},
\S \ref{SectionWeightSpectralSequences}).

The weight-monodromy conjecture claims the coincidence
of these two filtrations up to shift
(\cite{HodgeI}, \cite{WeilII}, \cite{Illusie},
\cite{Rapoport-Zink}).
Sometimes, the weight-monodromy conjecture
is also called {\it Deligne's conjecture on the purity
of monodromy filtration} in the literature.

\begin{conj}[Weight-monodromy conjecture]
\label{WMC}
$M_i V = W_{w+i} V$ for all $i$.
\end{conj}

In this paper, we prove Conjecture \ref{WMC}
for certain threefolds.
The main theorem of this paper is as follows.

\begin{thm}
\label{MainTheorem}
Let $X$ be a proper smooth variety of dimension 3 over $K$
which has a proper strictly semistable model $\X$ over $\O_K$.
Let $X_1,\ldots,X_m$ be the irreducible components
of the special fiber of $\X$.
Assume that the following conditions are satisfied.
\begin{enumerate}
\item $\X$ is projective over $\O_K$.
\item For all $i$, the Picard number $\rho((X_i)_{\overline{F}})$
is equal to the second $l$-adic Betti number
$\dim_{\Q_l} H^2((X_i)_{\overline{F}},\Q_l)$,
where $(X_i)_{\overline{F}} = X_i \otimes_F \overline{F}$.
\end{enumerate}
Then, Conjecture \ref{WMC} holds for $X$.
\end{thm}

Similarly, in Theorem \ref{MainTheoremPadic},
we also prove a $p$-adic analogue
by using the weight spectral sequence of Mokrane
for $p$-adic cohomology
(for details, see \S \ref{SectionPadic}).

For a proper smooth variety $Y$ over $\overline{F}$,
let $\Pic(Y)$ be the Picard group of $Y$,
and $\Pic^{\tau}(Y)$ (resp. $\Pic^0(Y)$)
the subgroup of $\Pic(Y)$ generated by divisors
which are numerically equivalent to 0
(resp. algebraically equivalent to 0).
The quotient $NS(Y) := \Pic(Y) / \Pic^0(Y)$ is called
the {\it N\'eron-Severi group} of $Y$.
It is known that $NS(Y)$ is a finitely generated
abelian group, $\Pic^{\tau}(Y) / \Pic^0(Y)$ is
a finite abelian group, and
$\Pic^0(Y)$ has a structure of an abelian variety
called the {\it Picard variety} of $Y$
(\cite{Fulton}, 19.3.1, \cite{FGA}, C-07 -- C-11).
$\rho(Y) := \text{rank}_{\Z}\,NS(Y)$ is called
the {\it Picard number} of $Y$.
Since algebraically equivalent cycles define
the same cohomology class,
we have the cycle map $NS(Y) \to H^2(Y,\Q_l(1))$.
If $x_1,\dots,x_r \in NS(Y) \otimes_{\Z} \Q$ are
linearly independent over $\Q$,
the images of $x_1,\dots,x_r$
by the cycle map is linearly independent over $\Q_l$
by an argument using Poincar\'e duality
(see also an argument in \cite{Kleiman}, Lemma 5-2).
Hence, the map
$NS(Y) \otimes_{\Z} \Q_l \to H^2(Y,\Q_l(1))$
is injective.
Therefore, the second condition of Theorem \ref{MainTheorem}
is equivalent to
say that the cycle map
$$
\begin{CD}
NS((X_i)_{\overline{F}}) \otimes_{\Z} \Q_l
    @>{\cong}>> H^2((X_i)_{\overline{F}},\Q_l(1)),
\end{CD}
$$
is an isomorphism for all $i$.

Conjecture \ref{WMC} was known to hold
if $X/K$ satisfies at least one of
the following conditions.
\begin{enumerate}
\item $X$ has a proper smooth model $\X$ over $\O_K$.
  In this case, Conjecture \ref{WMC} follows from
  the Weil conjecture for the special fiber
  $\X \otimes_{\O_K} F$ (\cite{WeilI},\cite{WeilII}).
\item $X$ is a curve or an abelian variety.
  In this case, Grothendieck proved Conjecture \ref{WMC}
  by the Picard-Lefschetz theorem and the theory of N\'eron models
  (\cite{SGA7-I},\ IX, for an arithmetic proof by Deligne,
  see also \cite{SGA7-I},\ I, 6, Appendice).
\item $X$ is of dimension $\leq 2$.
  If $X$ has a proper strictly semistable model $\X$ over $\O_K$,
  then Conjecture \ref{WMC} was proved by Rapoport-Zink
  by using the weight-spectral sequence
  (\S \ref{SectionWeightSpectralSequences}, \cite{Rapoport-Zink}).
  For general $X$, we may use de Jong's alteration (\cite{deJong}).
\item $K$ is of characteristic $p>0$.
  If $X/K$ comes from a family of varieties
  over a curve $f : \X \to C$ over a finite field
  by localization, Conjecture \ref{WMC} was proved
  by Deligne in \cite{WeilII}
  in his proof of the Weil conjecture.
  We can reduce the general case to Deligne's case
  by a specialization argument (\cite{Terasoma}, \cite{Ito}).
\item $K$ and $F$ are of characteristic $0$.
  In this case, Conjecture \ref{WMC} follows from
  a Hodge analogue over $\C$ by Lefschetz principle
  (Remark \ref{RemarkHodge}, \cite{Steenbrink}, \cite{MSaito1},
  see also \cite{SaitoZucker}, \cite{GuillenNavarroAznar},
  \cite{MSaito2}). In \cite{Ito}, another proof was given
  by modulo $p$ reduction.
\end{enumerate}

Our proof of Theorem \ref{MainTheorem} is based on a careful analysis
of the weight spectral sequence of Rapoport-Zink (\cite{Rapoport-Zink}).
It is inspired by the proof of a Hodge analogue over $\C$
by Steenbrink, M. Saito.
Instead of Hodge theory, we use the Hodge index theorem
for surfaces which was proved in any characteristic
by Segre and, independently, by Grothendieck
(see \cite{Kleiman}, 5, \cite{Grothendieck},
\cite{Fulton}, Example 15.2.4).

Note that, in dimension $\geq 3$ and in mixed characteristic,
Conjecture \ref{WMC} is still open in general.
In the followings, we give a way
to construct nontrivial examples for which the conditions of
Theorem \ref{MainTheorem} are satisfied
by using the theory of $p$-adic uniformization
by Mumford, Mustafin, Kurihara
(for details about $p$-adic uniformization, see
\cite{Mustafin}, \cite{Kurihara}).

\begin{ex}
Let $p$ be a prime number, $d \geq 3$ an integer,
$K$ a finite extension of the $p$-adic field $\Q_p$
with residue field $\F_q$,
$\widehat{\Omega}^d_K$ the Drinfeld upper half space
of dimension $d$ over $K$.
For a sufficiently small cocompact torsion free discrete subgroup
$\Gamma \subset \PGL_{d+1}(K)$, it is known that
the quotient $\Gamma \backslash \widehat{\Omega}^d_K$
can be algebraized to a projective smooth variety $X_{\Gamma}$ over $K$.
Moreover, $X_{\Gamma}$ has
a projective strictly semistable model $\X_{\Gamma}$ over $\O_K$
such that all irreducible components of
the special fiber are isomorphic to the blow-up of
$\P^d_{\F_q}$ along $\F_q$-rational linear subvarieties.
Therefore, if $d=3$, $X_{\Gamma}$ itself satisfies the conditions of
Theorem \ref{MainTheorem}.
More generally, by the hard Lefschetz theorem,
the conditions of Theorem \ref{MainTheorem} are satisfied
for the generic fiber of $\X_{\Gamma} \cap H_1 \cap \cdots \cap H_{d-3}$,
where $\X_{\Gamma} \cap H_1 \cap \cdots \cap H_{d-3}$
denotes the scheme obtained by taking general hyperplane
sections $(d-3)$-times
for a projective embedding $\X_{\Gamma} \subset \P^N_{\O_K}$.
\end{ex}

\begin{rem}
Recently, the author proved Conjecture \ref{WMC} for
$X_{\Gamma}$ for all $d$ (\cite{Ito2}).
The crucial point in \cite{Ito2} is to prove Grothendieck's
Hodge standard conjecture for varieties appearing
in the special fiber of $\X_{\Gamma}$.
\end{rem}

\vspace{0.2in}

\noindent
{\bf Acknowledgments.}
The author is grateful to Takeshi Saito and Kazuya Kato
for their advice and support.
He would like to thank
Laurent Clozel, Luc Illusie, Minhyong Kim, Barry Mazur,
Michael Rapoport, Richard Taylor, Teruyoshi Yoshida
for their interests in this work and encouragement,
Yoichi Mieda for reading manuscripts carefully
and pointing out errors in an earlier version,
and the referee for
invaluable comments and suggestions to make this paper
much more readable.
A part of this work was done during the author's stay
at Universit\'e de Paris-Sud in June 2001
and Korea Institute for Advanced Study in November 2001.
He would like to thank them for their cordial hospitality.
The author was supported by the Japan Society for the
Promotion of Science Research Fellowships for Young
Scientists.

\section{Monodromy filtration and weight filtration}
\label{SectionFiltration}

Let notation be as in \S \ref{SectionIntroduction}.
In this section, we recall the definition of
the monodromy filtration $M$ on $V$
and the weight filtration $W$ on $V$.

\subsection{Monodromy filtration}

Let $I_K$ be the inertia group of $K$,
which is the subgroup of $\Gal(\overline{K}/K)$ consisting of
elements acting trivially on $\overline{F}$.
Then there is an exact sequence
$$
\begin{CD}
1 @>>> I_K @>>> \Gal(\overline{K}/K) @>>> \Gal(\overline{F}/F) @>>> 1.
\end{CD}
$$
$\Gal(\overline{F}/F)$ acts on $I_K$ by conjugation
($\tau : \sigma \mapsto \tau \sigma \tau^{-1},
\ \tau \in \Gal(\overline{F}/F),\ \sigma \in I_K$).
The pro-$l$-part of $I_K$ is isomorphic to $\Z_l(1)$
as a $\Gal(\overline{F}/F)$-module by
$$ t_l : I_K \ni \sigma \mapsto
   \left( \frac{\sigma(\pi^{1/l^n})}{\pi^{1/l^n}} \right)_{\!\!n}
   \in \varprojlim \mu_{l^n} = \Z_l(1), $$
where $\pi$ is a uniformizer of $K$,
and $\mu_{l^n}$ is the group of $l^n$-th roots of unity.
It is known that $t_l$ is independent of
a choice of $\pi$ and its $l^n$-th root $\pi^{1/l^n}$ (\cite{Serre}).

By Grothendieck's monodromy theorem (\cite{SerreTate}, Appendix,
see also \cite{SGA7-I},\ I, Variante 1.3),
there exist $N,M \geq 1$ such that $(\rho(\sigma)^N - 1)^M = 0$
for all $\sigma \in I_K$.
Therefore, by replacing $K$ by a finite extension of it,
we may suppose that $I_K$ acts on $V$ through
$t_l : I_K \to \Z_l(1)$ and this action is unipotent.
Then there is a unique nilpotent map
$N : V(1) \to V$ of $\Gal(\overline{K}/K)$-representations
such that
$$ \rho(\sigma)
    = \exp(t_l(\sigma) N) \quad \text{for all} \quad \sigma \in I_K. $$
Here $N$ is a nilpotent map means that
there is $r \geq 1$ such that $N^r : V(r) \to V$ is zero,
where $V(r)$ denotes the $r$-th Tate twist of $V$.
$N$ is called the {\it monodromy operator} on $V$.

\begin{defn} (\cite{WeilII}, I, 1.7.2)
\label{DefMonodromyFiltration}
There exists a unique filtration $M$ on $V$
called the {\it monodromy filtration} characterized by
the following properties.
\begin{enumerate}
\item $M$ is an increasing filtration
  $\cdots \subset M_{i-1} V \subset M_i V \subset M_{i+1} V \subset \cdots$
  of $\Gal(\overline{K}/K)$-representations such that
  $M_i V = 0$ for sufficiently small $i$ and
  $M_i V = V$ for sufficiently large $i$.
\item $N (M_i V(1)) \subset M_{i-2} V$ for all $i$.
\item By the second condition,
  we can define $N : \Gr^M_{i} V (1) \to \Gr^M_{i-2} V$,
  where $\Gr^M_{i} V = M_{i} V / M_{i-1} V$.
  Then, for each $r \geq 0$, $N^r : \Gr^M_{r} V (r) \to \Gr^M_{-r} V$
  is an isomorphism.
\end{enumerate}
\end{defn}

In the above definition,
we have replaced $K$ by a finite extension of it.
We can easily see that
$M$ is stable under the action of $\Gal(\overline{K}/K)$
for the original $K$.
Therefore, we can define the monodromy filtration $M$
as a filtration of $\Gal(\overline{K}/K)$-representations
without replacing $K$ by a finite extension of it.

\subsection{Weight filtration}

Here we give two definitions of the weight filtration $W$ on $V$.
One definition makes sense if $F$ is a finite field
(Definition \ref{DefWeightFiltration}).
The other one makes sense if $X$ has
a proper strictly semistable model over $\O_K$
(Definition \ref{DefWeightFiltration2}).
Note that if both definitions
make sense, they define the same filtration
(Remark \ref{EquivalenceTwoDefinitions},
Remark \ref{PropertiesWeightSpectralSequence}, 2).

Firstly, assume that $F$ is a finite field $\F_q$
with $q$ elements.
Let $\Fr_q \in \Gal(\overline{\F}_q/\F_q)$ be
the geometric Frobenius element.
Namely, $\Fr_q$ is the inverse of the $q$-th power map
($\overline{\F}_q \ni x \mapsto x^{q} \in \overline{\F}_q$).

\begin{defn} (\cite{HodgeI}, \cite{WeilII}, I, 1.7.5)
\label{DefWeightFiltration}
There exists a unique filtration $W$ 
called the {\it weight filtration} on $V$ characterized by
the following properties
(for existence, see Remark \ref{PropertiesWeightSpectralSequence}, 2).
\begin{enumerate}
\item $W$ is an increasing filtration
  $\cdots \subset W_{i-1} V \subset W_i V \subset W_{i+1} V \subset \cdots$
  of $\Gal(\overline{K}/K)$-representations such that
  $W_i V = 0$ for sufficiently small $i$ and
  $W_i V = V$ for sufficiently large $i$.
\item The action of $I_K$ on $\Gr^W_{i} V$ factors through a finite quotient.
\item Let $\widetilde{\Fr}_q$ be a lift of $\Fr_q$ in $\Gal(\overline{K}/K)$.
  Then all of the eigenvalues of the action of $\widetilde{\Fr}_q$
  on $\Gr^W_{i} V = W_i V / W_{i-1} V$ are algebraic integers
  with the property that the complex absolute values of their complex
  conjugates are $q^{i/2}$.
  Note that this condition doesn't depend on
  a choice of $\widetilde{\Fr}_q$ by the second condition.
\end{enumerate}
\end{defn}

If $X$ has a proper strictly semistable model over $\O_K$
(Definition \ref{DefSemistableModel}),
we define the weight filtration $W$ on $V$
by using the weight spectral sequence of Rapoport-Zink
(for details, see \S \ref{SectionWeightSpectralSequences}).

\begin{defn}
\label{DefWeightFiltration2}
The {\it weight filtration} $W$ on $V$ is a filtration
defined by the weight spectral sequence of Rapoport-Zink
in \S \ref{SectionWeightSpectralSequences}.
\end{defn}

\begin{rem}
\label{EquivalenceTwoDefinitions}
By the Weil conjecture, if $F$ is a finite field
and $X$ has a proper strictly semistable model over $\O_K$,
these two definitions give the same filtration
on $V$ (Remark \ref{PropertiesWeightSpectralSequence}, 2).
\end{rem}

\section{Weight spectral sequences}
\label{SectionWeightSpectralSequences}

Let notation be as in \S \ref{SectionIntroduction}.
Let $X$ be a proper smooth variety over $K$ of dimension $n$.

\begin{defn}
\label{DefSemistableModel}
A regular scheme $\X$ which is proper and flat over $\O_K$
is called a {\em proper strictly semistable model}
of $X$ over $\O_K$ if the generic fiber
$\X_K := \X \otimes_{\O_K} K$ is isomorphic to $X$
and the special fiber
$\X_{F} := \X \otimes_{\O_K} F$
is a divisor of $\X$ with simple normal crossings.
\end{defn}

We recall the weight spectral sequence of
Rapoport-Zink (\cite{Rapoport-Zink}).
Assume that $X$ has a proper strictly semistable model $\X$ over $\O_K$.
Let $X_1,\ldots,X_m$ be the irreducible components of
the special fiber of $\X$, and
$$ X^{(j)} := \coprod_{1 \leq i_1 < \cdots < i_j \leq m}
     X_{i_1} \cap \cdots \cap X_{i_j}. $$
Then $X^{(j)}$ is a proper smooth variety of dimension $n-j+1$ over $F$.
The {\it weight spectral sequence of Rapoport-Zink}
is as follows.
\begin{center}
\begin{tabular}{l}
$\displaystyle E_1^{-r,\,w+r} = \bigoplus_{k \geq \max\{0,-r\}}
      H^{w-r-2k}(X^{(2k+r+1)}_{\overline{F}},\, \Q_l(-r-k))$ \\
\hspace*{3.2in} $\displaystyle \Longrightarrow \quad
      H^w(X_{\overline{K}},\Q_l)$
\end{tabular}
\end{center}
This spectral sequence is $\Gal(\overline{K}/K)$-equivariant.
The map $d_1^{i,\,j} : E_1^{i,\,j} \to E_1^{i+1,\,j}$ can be described
by restriction morphisms and Gysin morphisms explicitly
(see \cite{Rapoport-Zink}, 2.10 for details).
Note that $I_K$ acts on $E_1^{i,\,j}$ trivially and
$\Gal(\overline{F}/F)$ acts on them.

\begin{rem}
\label{PropertiesWeightSpectralSequence}
We recall some properties of the weight spectral sequence.
\begin{enumerate}
\item The action of the monodromy operator $N$
  on $H^w(X_{\overline{K}},\Q_l)$
  is induced by a natural map $N : E_1^{i,\,j}(1) \to E_1^{i+2,\,j-2}$
  satisfying
    $$ \begin{CD}
    N^r : E_1^{-r,\,w+r}(r) @>{\cong}>> E_1^{r,\,w-r}.
    \end{CD} $$
  for all $r,w$.
  We can describe $N : E_1^{i,\,j}(1) \to E_1^{i+2,\,j-2}$
  explicitly (\cite{Rapoport-Zink}, 2.10).
\item If $F$ is a finite field $\F_q$ with $q$ elements,
    $ H^{w-r-2k}(X^{(2k+r+1)}_{\overline{F}},\, \Q_l(-r-k)) $
  has {\it weight} $(w-r-2k) -2(-r-k) = w+r$
  by the Weil conjecture (\cite{WeilI}, \cite{WeilII}).
  Namely, all of the eigenvalues of the action of $\Fr_q$
  are algebraic integers with the property that
  the complex absolute values of their complex
  conjugates are $q^{(w+r)/2}$.
  Therefore, $E_1^{i,\,j}$ has weight $j$, and
  the filtration on $H^w(X_{\overline{K}},\Q_l)$
  induced by the weight spectral sequence is
  the weight filtration $W$ 
  in Definition \ref{DefWeightFiltration}.
  Note that this proves the existence of $W$ in
  Definition \ref{DefWeightFiltration}
  (for general $X$, we may use
  de Jong's alteration to reduce to the semistable case
  (\cite{deJong})).
\item The weight spectral sequence always degenerates at $E_2$.
  If $F$ is a finite field, this is a consequence of
  the Weil conjecture.
  Since $E_1^{i,\,j}$ has weight $j$ as above,
  $d_r^{i,\,j} : E_r^{i,\,j} \to E_r^{i+r,\,j-r+1}$
  is a map between $\Gal(\overline{F}/F)$-representations
  with different weights for $r \geq 2$.
  Hence it must be zero.
  For general $F$, Nakayama proved the $E_2$-degeneracy
  by reducing to the above case by log geometry
  (\cite{CNakayama},
  for the equal characteristic case, see also \cite{Ito}).
\end{enumerate}
\end{rem}

Therefore,
if $X$ has a proper strictly semistable model over $\O_K$,
the weight-monodromy conjecture (Conjecture \ref{WMC})
is equivalent to the following conjecture
on the weight spectral sequence.

\begin{conj}[\cite{Rapoport-Zink}, \cite{Illusie}]
\label{WMC_semistable}
Let $X$ be a proper smooth variety over $K$
which has a proper strictly semistable model $\X$ over $\O_K$.
Let $E_1^{-r,\,w+r} \Rightarrow H^w(X_{\overline{K}},\Q_l)$
be the weight spectral sequence of Rapoport-Zink.
Then $N^r$ induces an isomorphism
$$
\begin{CD}
N^r : E_2^{-r,\,w+r}(r) @>{\cong}>> E_2^{r,\,w-r}
\end{CD}
$$
on $E_2$-terms for all $r,w$.
\end{conj}

\begin{rem}
\label{RemarkHodge}
The weight spectral sequence was originally constructed
by Steenbrink over $\C$ (\cite{Steenbrink}).
In this case, we can formulate a Hodge analogue of
Conjecture \ref{WMC_semistable}.
This was proved by Steenbrink, M. Saito,
and Deligne (unpublished)
by using polarized Hodge structures
(\cite{Steenbrink}, 5.10, \cite{MSaito1}, 4.2.5,
see also \cite{SaitoZucker}, \cite{GuillenNavarroAznar},
\cite{MSaito2}).
If $F$ is of characteristic 0,
we can prove Conjecture \ref{WMC_semistable}
by Lefschetz principle.
However, in mixed characteristic,
we can't directly follow the argument over $\C$
because we don't have a good analogue of polarized
Hodge structures for \'etale cohomology.
Nevertheless, for divisors on algebraic surfaces,
we have a good analogue, namely,
the Hodge index theorem.
This is a crucial observation in the proof
of Theorem \ref{MainTheorem} in this paper
(see \cite{MSaito2} for an argument assuming
Grothendieck's standard conjectures).
\end{rem}

\section{Proof of Theorem \ref{MainTheorem}}

Let notation be as in \S \ref{SectionIntroduction}.
Let $X$ be a proper smooth variety of dimension 3 over $K$
which has a proper strictly semistable model $\X$ over $\O_K$.
Let $X_1,\ldots,X_m$ be the irreducible components of
the special fiber of $\X$.
Let $E_1^{-r,\,w+r} \Rightarrow H^w(X_{\overline{K}},\Q_l)$
be the weight spectral sequence of Rapoport-Zink
as in \S \ref{SectionWeightSpectralSequences}.
To prove Theorem \ref{MainTheorem},
it is enough to prove the following proposition
(see Conjecture \ref{WMC_semistable}).

\begin{prop}
\label{WMC3fold}
Let notation be as above.
Assume that the following conditions are satisfied.
\begin{enumerate}
\item $\X$ is projective over $\O_K$.
\item $\rho((X_i)_{\overline{F}}) =
  \dim_{\Q_l} H^2((X_i)_{\overline{F}},\Q_l)$ for all $i$.
\end{enumerate}
Then $N^r$ induces an isomorphism
$$
\begin{CD}
N^r : E_2^{-r,\,w+r}(r) @>{\cong}>> E_2^{r,\,w-r}
\end{CD}
$$
on $E_2$-terms for all $r,w$.
Namely, Conjecture \ref{WMC_semistable} holds for $\X$.
Hence Conjecture \ref{WMC} also holds for $X$.
\end{prop}

\begin{proof}
First of all, we simplify and fix some notation.
Let
$$
\begin{CD}
\displaystyle
d_1^{i,j} = \sum_{k \geq \max\{0,-r\}} (\rho + \tau) : E_1^{i,j}
  @>>> E_1^{i+1,j}
\end{CD}
$$
denote the differential on $E_1$-terms, where
$$
\begin{CD}
\rho : H^s(X^{(t)}_{\overline{F}},\Q_l) @>>>
  H^s(X^{(t+1)}_{\overline{F}},\Q_l)
\end{CD}
$$
be a linear combination of restriction morphisms for some $s,t$
($\rho$ is $(-1)^{r+k} \theta$ in \cite{Rapoport-Zink}, 2.10),
and
$$
\begin{CD}
\tau : H^s(X^{(t)}_{\overline{F}},\Q_l) @>>>
  H^{s+2}(X^{(t-1)}_{\overline{F}},\Q_l(1))
\end{CD}
$$
be a linear combination of Gysin morphisms for some $s,t$
($\tau$ is $(-1)^k d_1'$ in \cite{Rapoport-Zink}, 2.10).
We use the same letters $\rho, \tau$ for different $i,j,k$,
by abuse of notation.
We write
$$ A = X^{(1)}_{\overline{F}}, \qquad
   B = X^{(2)}_{\overline{F}}, \qquad
   C = X^{(3)}_{\overline{F}}, \qquad
   D = X^{(4)}_{\overline{F}}. $$
Note that $A,B,C,D$ are disjoint unions of proper smooth
varieties of dimension $3,2,1,0$, respectively.
We omit the coefficients and the Tate twists of
\'etale cohomology groups.

Then we have the following sequence of restriction morphisms
$$
\begin{CD}
H^s(A) @>{\rho}>> H^s(B) @>{\rho}>> H^s(C) @>{\rho}>> H^s(D)
\end{CD}
$$
satisfying $\rho \circ \rho = 0$, and
the following sequence of Gysin morphisms
$$
\begin{CD}
H^s(D) @>{\tau}>> H^{s+2}(C) @>{\tau}>> H^{s+4}(B) @>{\tau}>> H^{s+6}(A)
\end{CD}
$$
satisfying $\tau \circ \tau = 0$.
Moreover, from the explicit description of $\rho, \tau$
in \cite{Rapoport-Zink}, 2.10,
we know that $\rho, \tau$ are dual to each other with respect to
the sum of the cup product pairings and
satisfy $\tau \circ \rho + \rho \circ \tau = 0$.

In the notation above,
the $E_1$-terms $E_1^{i,j}$ of the weight spectral sequence
are as follows.
\begin{center}
{
\noindent
\hspace*{-0.2in}
\begin{tabular}{llllllllllllll}
{\scriptsize $6$}
 & $H^0(D)$ & \bb $\stackrel{\tau}{\rightarrow}$ &
  \bb $H^2(C)$ & \bb $\stackrel{\tau}{\rightarrow}$ &
  \bb $H^4(B)$ & \bb $\stackrel{\tau}{\rightarrow}$ &
  \bb $H^6(A)$ \vspace{4mm} \\
{\scriptsize $5$}
 & & & \bb $H^1(C)$ & \bb $\stackrel{\tau}{\rightarrow}$ &
  \bb $H^3(B)$ & \bb $\stackrel{\tau}{\rightarrow}$ &
  \bb $H^5(A)$ \vspace{3mm} \\
{\scriptsize $4$}
 & & & \bb $H^0(C)$ & \bb $\stackrel{\tau+\rho}{\rightarrow}$ &
   \hspace*{-6mm}
    \begin{tabular}{ll} $H^2(B)$ \\ $\oplus H^0(D)$ \end{tabular}
   \hspace*{-4mm}
  & \bb $\stackrel{(\tau+\rho, \tau)}{\rightarrow}$ &
   \hspace*{-6mm}
    \begin{tabular}{ll} $H^4(A)$ \\ $\oplus H^2(C)$ \end{tabular}
   \hspace*{-4mm}
  & \bb $\stackrel{(\rho,\tau)}{\rightarrow}$ &
  \bb $H^4(B)$ \vspace{3mm} \\
{\scriptsize $3$}
 & & & & & \bb $H^1(B)$ & \bb $\stackrel{\tau+\rho}{\rightarrow}$ &
   \hspace*{-6mm}
    \begin{tabular}{ll} $H^3(A)$ \\ $\oplus H^1(C)$ \end{tabular}
   \hspace*{-4mm}
  & \bb $\stackrel{(\rho,\tau)}{\rightarrow}$ &
  \bb $H^3(B)$ \vspace{3mm} \\ 
{\scriptsize $2$}
 & & & & & \bb $H^0(B)$ & \bb $\stackrel{\tau+\rho}{\rightarrow}$ &
   \hspace*{-6mm}
    \begin{tabular}{ll} $H^2(A)$ \\ $\oplus H^0(C)$ \end{tabular}
   \hspace*{-4mm}
  & \bb $\stackrel{(\rho,\tau+\rho)}{\rightarrow}$ &
   \hspace*{-6mm}
    \begin{tabular}{ll} $H^2(B)$ \\ $\oplus H^0(D)$ \end{tabular}
   \hspace*{-4mm}
  & \bb $\stackrel{(\rho,\tau)}{\rightarrow}$ &
  \bb $H^2(C)$ \vspace{3mm} \\
{\scriptsize $1$}
 & & & & & & & \bb $H^1(A)$ & \bb $\stackrel{\rho}{\rightarrow}$ &
  \bb $H^1(B)$ & \bb $\stackrel{\rho}{\rightarrow}$ &
  \bb $H^1(C)$ \vspace{4mm} \\
{\scriptsize $0$}
 & & & & & & & \bb $H^0(A)$ & \bb $\stackrel{\rho}{\rightarrow}$ &
  \bb $H^0(B)$ & \bb $\stackrel{\rho}{\rightarrow}$ &
  \bb $H^0(C)$ & \bb $\stackrel{\rho}{\rightarrow}$ &
  \bb $H^0(D)$ \\
{\scriptsize $j/i$}
 & {\scriptsize $-3$} & & {\scriptsize $-2$} & &
   {\scriptsize $-1$} & & {\scriptsize $0$} & &
   {\scriptsize $1$}  & & {\scriptsize $2$} & &
   {\scriptsize $3$}
\end{tabular}
}
\end{center}

We have to check that $N^r$ induces an isomorphism
$N^r : E_2^{-r,\,w+r} \stackrel{\cong}{\to} E_2^{r,\,w-r}$
on $E_2$-terms for all $r,w$.

For most $r,w$, we can use the same argument as in
\cite{Rapoport-Zink}, 2.13.
Namely, for $H^0$ and the Poincar\'e dual of $H^0$,
the cohomology groups have combinatorial structures.
Hence we can apply a positivity argument
to the underlying $\Q$-structures of them.
For $H^1$ and the Poincar\'e dual of $H^1$,
we can use polarizations of Picard varieties
(for details of the proof,
see \cite{Rapoport-Zink}, 2.13).

The only part to which we can't apply
an argument in \cite{Rapoport-Zink} is the following part.
$$
\xymatrix{
H^0(C) \ar[r]^-{\tau+\rho} &
   \text{\begin{tabular}{ll} $H^2(B)$ \\ $\oplus H^0(D)$ \end{tabular}}
   \ar[r]^-{(\tau+\rho, \tau)} \ar[ddrr]^{N} &
   \text{\begin{tabular}{ll} $H^4(A)$ \\ $\oplus H^2(C)$ \end{tabular}} \\
\\
& & \text{\begin{tabular}{ll} $H^2(A)$ \\ $\oplus H^0(C)$ \end{tabular}}
  \ar[r]^-{(\rho,\tau+\rho)} &
  \text{\begin{tabular}{ll} $H^2(B)$ \\ $\oplus H^0(D)$ \end{tabular}}
  \ar[r]^-{(\rho,\tau)} &
  H^2(C)
}
$$
Here $N$ is the identity map.
Note that these two rows are dual to each other with respect to
the sum of the cup product pairings
\begin{align*}
H^2(A) \times H^4(A) &\to \Q_l, \qquad
H^2(B) \times H^2(B) \to \Q_l \\
H^0(C) \times H^2(C) &\to \Q_l, \qquad
H^0(D) \times H^0(D) \to \Q_l.
\end{align*}
By definition, the cohomology of the first row is $E_2^{-1,4}$,
and the cohomology of the second row is $E_2^{1,2}$.

Let us consider the following simple linear algebra lemma.

\begin{lem}
\label{LemmaLinearAlgebra}
Let
$$
\begin{CD}
V_1 @>{f}>> V_2 @>{g}>> V_3
\end{CD}
$$
be a complex of finite dimensional vector spaces over a field $k$,
and
$$
\begin{CD}
V_3^{\ast} @>{g^{\ast}}>> V_2^{\ast} @>{f^{\ast}}>> V_1^{\ast}
\end{CD}
$$
the dual of it.
Fix a nondegenerate pairing $\langle, \rangle$ on $V_2$,
and identify $V_2$ and $V_2^{\ast}$ by $\langle, \rangle$.
Assume that $\Im f \subset \Im g^{\ast}$.
Then, the identification $V_2 = V_2^{\ast}$
induces a map on the cohomology groups
$$
\begin{CD}
\Ker g / \Im f @>>> \Ker f^{\ast} / \Im g^{\ast}.
\end{CD}
$$
This is an isomorphism
if and only if $\Ker g \cap \Im g^{\ast} \subset \Im f$.
\end{lem}

\begin{proof}[Proof of Lemma \ref{LemmaLinearAlgebra}]
Since
\begin{align*}
  \dim_{k} (\Ker g / \Im f)
     &= (\dim_{k} V_2 - \text{rank}\,g) - \text{rank}\,f \\
  \dim_{k} (\Ker f^{\ast} / \Im g^{\ast})
     &= (\dim_{k} V_2^{\ast} - \text{rank}\,f^{\ast})
           - \text{rank}\,g^{\ast},
\end{align*}
we see that $\Ker g / \Im f$ and $\Ker f^{\ast} / \Im g^{\ast}$
have the same dimensions.
Therefore,
$\Ker g / \Im f \to \Ker f^{\ast} / \Im g^{\ast}$
is an isomorphism if and only if it is injective,
if and only if $\Ker g \cap \Im g^{\ast} \subset \Im f$.
\end{proof}

To prove Proposition \ref{WMC3fold},
we need the following important lemma.

\begin{lem}
\label{KeyLemma}
In the following sequence of $\tau, \rho, \tau$
$$
\begin{CD}
H^0(B) @>{\tau}>> H^{2}(A) @>{\rho}>> H^{2}(B) @>{\tau}>> H^{4}(A),
\end{CD}
$$
we have
$\Ker \tau \cap \Im \rho = \Im(\rho \circ \tau)$
in $H^2(B)$.
\end{lem}

Note that since 
$\tau \circ \rho \circ \tau = - \rho \circ \tau \circ \tau = 0$,
$\Ker \tau \cap \Im \rho \supset \Im(\rho \circ \tau)$.
However the opposite inclusion $\subset$ is highly nontrivial.

We postpone the proof of Lemma \ref{KeyLemma}
and finish the proof of Proposition \ref{WMC3fold}.
According to Lemma \ref{LemmaLinearAlgebra},
it is enough to show that the intersection
\begin{align*}
&\Ker \big( H^2(B) \oplus H^0(D)
    \stackrel{(\tau+\rho, \tau)}{\longrightarrow}
    H^4(A) \oplus H^2(C) \big) \\
&\hspace*{1in} \cap \ \Im \big( H^2(A) \oplus H^0(C)
    \stackrel{(\rho, \tau+\rho)}{\longrightarrow}
    H^2(B) \oplus H^0(D) \big)
\end{align*}
is contained in
$\Im \big( H^0(C) \stackrel{\tau+\rho}{\longrightarrow}
    H^2(B) \oplus H^0(D) \big)$.
Take an element
$(a,c) \in H^2(A) \oplus H^0(C)$
such that the image $(\rho(a) + \tau(c),\,\rho(c))$
lies in the above intersection.
Namely, we have
$$ \tau(\rho(a) + \tau(c)) = 0, \qquad
   \rho(\rho(a) + \tau(c)) + \tau \circ \rho(c) = 0. $$
Since $\tau \circ \tau = 0$, we have
$\rho(a) \in \Ker \tau \cap \Im \rho \subset H^2(B)$.
Hence by Lemma \ref{KeyLemma},
we can write $\rho(a) = \rho \circ \tau(b)$ for some $b \in H^0(B)$.
We put $c' = - \rho(b) + c \in H^0(C)$.
Then we have
$$
   \tau(c') = \tau(- \rho(b) + c) = \rho(a) + \tau(c), \qquad
   \rho(c') = \rho(- \rho(b) + c) = \rho(c)
$$
since $- \tau \circ \rho = \rho \circ \tau,\ \rho \circ \rho = 0$.
Hence we have
$$ (\rho(a)+\tau(c),\,\rho(c))
      \in \Im \big( H^0(C) \stackrel{\tau+\rho}{\longrightarrow}
          H^2(B) \oplus H^0(D) \big). $$
and Proposition \ref{WMC3fold} is proved
except for Lemma \ref{KeyLemma}.
\end{proof}

Finally we shall prove Lemma \ref{KeyLemma}.

\begin{proof}[Proof of Lemma \ref{KeyLemma}]
By the first condition of Proposition \ref{WMC3fold},
let
$$ L : H^s(A) \to H^{s+2}(A), \qquad
   L : H^s(B) \to H^{s+2}(B) $$
denote the cup products with the cohomology class of
an ample line bundle on the special fiber of $\X$.
$L$ commutes with $\rho, \tau$.
By the hard Lefschetz theorem
proved by Deligne in any characteristic
(\cite{WeilII}, 4.1.1), the maps
$$ L^i : H^{3-i}(A) \to H^{3+i}(A), \qquad
   L^i : H^{2-i}(A) \to H^{2+i}(B) $$
are isomorphisms.
We define the primitive parts as follows.
\begin{align*}
   P^0(A) &:= H^0(A), \quad
   P^2(A) := \Ker(L^2:H^2(A) \to H^6(A)), \\
   P^0(B) &:= H^0(B), \quad
   P^2(B) := \Ker(L:H^2(B) \to H^4(B)).
\end{align*}
To distinguish different $\tau$ and $\rho$,
we use the notation as in the following diagram.
$$
\xymatrix@=0pt@C=-60pt@R=0pt{
**[l] L^3 P^0(A) = H^6(A) \\
          & **[r] \ar[lu]_{\tau_4} H^4(B) = L^2 P^0(B) \\
**[l] L P^2(A) \oplus L^2 P^0(A) = H^4(A) \ar[ru]^{\rho_4} \\
          & **[r] \ar[lu]_{\tau_2} H^2(B) = L P^0(B) \oplus P^2(B)\\
**[l] P^2(A) \oplus L P^0(A) = H^2(A) \ar[ru]^{\rho_2} \\
          & **[r] \ar[lu]_{\tau_0} H^0(B) = P^0(B) \\
**[l] P^0(A) = H^0(A) \ar[ru]^{\rho_0} \\
}
$$
For $i=0,2,4$, we define
$\Im^{\!0} \rho_i \subset \Im \rho_i$, $\Im^{\!0} \tau_i \subset \Im \tau_i$
as follows.
$$
\begin{array}{lll}
\Im^{\!0} \rho_0 &:=& \Im \rho_0 \\
\Im^{\!0} \rho_2 &:=& (\rho_2(L P^0(A)) \cap L P^0(B)) \\
   & & \qquad \oplus (\Im \rho_2 \cap P^2(B)) \\
   &=& L \Im \rho_0 \oplus (\Im \rho_2 \cap P^2(B)) \\
\Im^{\!0} \rho_4 &:=& \rho_4(L^2 P^0(A)) \\
   &=& L^2 \Im \rho_0
\end{array} \quad
\begin{array}{lll}
\Im^{\!0} \tau_0 &:=& \Im \tau_0 \cap P^2(A) \\
\Im^{\!0} \tau_2 &:=& \tau_2(L P^0(B)) \cap L P^2(A) \\
   &=& L \Im^{\!0} \tau_0 \\
\Im^{\!0} \tau_4 &:=& 0 \\
\end{array}
$$
For $i=0,2,4$, we define
$$ \Im^{\!1} \rho_i := \Im \rho_i / \Im^{\!0} \rho_i,
\qquad \Im^{\!1} \tau_i := \Im \tau_i / \Im^{\!0} \tau_i. $$

\begin{claim}
\label{LemmaClaim1}
$L$ and $L^2$ induce the following isomorphisms.
$$ L^2 : \Im^{\!0} \rho_0 \cong \Im^{\!0} \rho_4, \qquad
   L : \Im^{\!0} \tau_0 \cong \Im^{\!0} \tau_2, $$
$$ L : \Im^{\!1} \rho_2 \cong \Im^{\!1} \rho_4, \qquad
   L^2 : \Im^{\!1} \tau_0 \cong \Im^{\!1} \tau_4. $$
\end{claim}

\begin{proof}[Proof of Claim \ref{LemmaClaim1}]
The claim for $\Im^{\!0}$ trivially follows from the definition
of $\Im^{\!0}$.
For $\Im^{\!1}$, the surjectivity follows from
the hard Lefschetz theorem.
And the injectivity can be checked directly
from the definition of $\Im^{\!0}$.
\end{proof}

\begin{claim}
\label{LemmaClaim2}
We have the following equality of dimensions of $\Im^{\!0}$
and $\Im^{\!1}$.
\begin{align*}
\dim_{\Q_l} \Im^{\!0} \rho_i &= \dim_{\Q_l} \Im^{\!1} \tau_i \qquad
   \text{\rm for}\ i=0,2,4 \\
\dim_{\Q_l} \Im^{\!0} \tau_i &= \dim_{\Q_l} \Im^{\!1} \rho_{i+2} \quad
   \text{\rm for}\ i=0,2
\end{align*}
\end{claim}

\begin{proof}[Proof of Claim \ref{LemmaClaim2}]
Since $\rho_0$ and $\tau_4$ are dual to each other,
$\dim_{\Q_l} \Im \rho_0 = \dim_{\Q_l} \Im \tau_4$.
Similarly,
$\dim_{\Q_l} \Im \rho_2 = \dim_{\Q_l} \Im \tau_2$,
$\dim_{\Q_l} \Im \rho_4 = \dim_{\Q_l} \Im \tau_0$.
On the other hand, since
$\dim_{\Q_l} \Im \rho_i
  = \dim_{\Q_l} \Im^{\!0} \rho_i + \dim_{\Q_l} \Im^{\!1} \rho_i$
and
$\dim_{\Q_l} \Im \tau_i
  = \dim_{\Q_l} \Im^{\!0} \tau_i + \dim_{\Q_l} \Im^{\!1} \tau_i$
for $i=0,2,4$,
we have
\begin{align*}
\dim_{\Q_l} \Im^{\!0} \rho_0 + \dim_{\Q_l} \Im^{\!1} \rho_0
  &= \dim_{\Q_l} \Im^{\!0} \tau_4 + \dim_{\Q_l} \Im^{\!1} \tau_4 \\
\dim_{\Q_l} \Im^{\!0} \rho_2 + \dim_{\Q_l} \Im^{\!1} \rho_2
  &= \dim_{\Q_l} \Im^{\!0} \tau_2 + \dim_{\Q_l} \Im^{\!1} \tau_2 \\
\dim_{\Q_l} \Im^{\!0} \rho_4 + \dim_{\Q_l} \Im^{\!1} \rho_4
  &= \dim_{\Q_l} \Im^{\!0} \tau_0 + \dim_{\Q_l} \Im^{\!1} \tau_0.
\end{align*}
Then Claim \ref{LemmaClaim2} follows from
Claim \ref{LemmaClaim1} and
$\dim_{\Q_l} \Im^{\!1} \rho_0 = \dim_{\Q_l} \Im^{\!0} \tau_4 = 0$.
\end{proof}

We define a pairing on $P^2(A)$ by
$P^2(A) \ni x,y \mapsto L \cup x \cup y \in \Q_l$,
where $L \cup x \cup y$ denotes the sum of
the cup product pairings.
We call this pairing {\it the sum of the Lefschetz pairings}
on $P^2(A)$.

\begin{claim}
\label{LemmaClaim3}
The restriction of the sum of the cup product pairings
on $H^2(B)$ to $\Im^{\!0} \rho_2$ is nondegenerate.
Moreover, the restriction of the sum of the Lefschetz pairings
on $P^2(A)$ to $\Im^{\!0} \tau_0$ is also nondegenerate.
\end{claim}

\begin{proof}[Proof of Claim \ref{LemmaClaim3}]
By the second condition of Proposition \ref{WMC3fold},
the cycle map
$$
\begin{CD}
NS(A) \otimes_{\Z} \Q_l @>{\cong}>> H^2(A)
\end{CD}
$$
is an isomorphism, where $NS(A)$ is the N\'eron-Severi group
of $A$.
We consider the following commutative diagram
$$
\begin{CD}
NS(A) \otimes_{\Z} \Q @>{\rho'}>> NS(B) \otimes_{\Z} \Q \\
@VVV @VVV \\
H^2(A) @>{\rho_2}>> H^2(B),
\end{CD}
$$
where $\rho'$ is defined by the same way as $\rho_2$.
We can define an analogue of the primitive decomposition
at the level of the N\'eron-Severi groups
$$ NS(A) \otimes_{\Z} \Q = N^0(A) \oplus N^2(A),\quad
   NS(B) \otimes_{\Z} \Q = N^0(B) \oplus N^2(B) $$
such that
\begin{align*}
 N^0(A) \otimes_{\Q} \Q_l &= L P^0(A), \quad
 N^2(A) \otimes_{\Q} \Q_l = P^2(A), \\
 N^0(B) \otimes_{\Q} \Q_l &= L P^0(B), \quad
 N^2(B) \otimes_{\Q} \Q_l \subset P^2(B).
\end{align*}
We define
$$ \Im^{\!0} \rho' := (\rho'(N^0(A)) \cap N^0(B))
       \oplus (\Im \rho' \cap N^2(B)). $$
Then we have
$(\Im^{\!0} \rho') \otimes_{\Q} \Q_l = \Im^{\!0} \rho_2$.
Namely, $\Im^{\!0} \rho_2$ has a $\Q$-structure
coming from the N\'eron-Severi groups.
Since the restriction of the cup product pairings on $H^2(B)$
to $NS(B) \otimes_{\Z} \Q$ is the intersection product
pairings, we have only to show
that the restriction of the sum of the intersection product
pairings on $NS(B) \otimes_{\Z} \Q$ to
$\Im^{\!0} \rho'$ is nondegenerate.
To prove the nondegeneracy,
take an element $x \in \Im^{\!0} \rho'$ such that
$x \cdot y = 0$ for all $y \in \Im^{\!0} \rho'$,
where $x \cdot y$ denotes
the sum of the intersection product pairings.
We write
$$ x = a + b,\qquad a \in \rho'(N^0(A)) \cap N^0(B),\quad
      b \in \Im \rho' \cap N^2(B) $$
If $a \neq 0$, we have $x \cdot a = a \cdot a > 0$
since the self intersection of an ample divisor is
positive.
Hence we have $a = 0$.
If $b \neq 0$, we have $x \cdot b = b \cdot b < 0$
by the Hodge index theorem
(see \cite{Kleiman}, 5, \cite{Grothendieck},
\cite{Fulton}, Example 15.2.4).
Therefore we have $b = 0$, hence $x = 0$.
This proves the first assertion.

The proof of the second assertion is similar.
By applying the Hodge index theorem
to a hyperplane section of
each connected component of $A$,
we see that the restriction of the sum of
the Lefschetz pairings to $N^2(A)$
is negative definite.
By the same way as above, we see that
$\Im^{\!0} \tau_0 \subset P^2(A)$ has a $\Q$-structure.
Namely, there is a $\Q$-subspace $V \subset N^2(A)$
such that $\Im^{\!0} \tau_0 = V \otimes_{\Q} \Q_l$.
The restriction of the sum of
the Lefschetz pairings to $V$ is nondegenerate
because it is negative definite.
Therefore, the restriction of the sum of
the Lefschetz pairings to $\Im^{\!0} \tau_0$
is also nondegenerate.
\end{proof}

\begin{claim}
\label{LemmaClaim4}
The composition of the following maps is an isomorphism.
$$
\begin{CD}
\Im^{\!0} \tau_0 @>{\rho_2}>> \Im \rho_2
       @>>> \Im^{\!1} \rho_2 = \Im \rho_2 / \Im^{\!0} \rho_2
\end{CD}
$$
\end{claim}

\begin{proof}[Proof of Claim \ref{LemmaClaim4}]
By Claim \ref{LemmaClaim2},
$\dim_{\Q_l} \Im^{\!0} \tau_0 = \dim_{\Q_l} \Im^{\!1} \rho_2$.
Therefore, we have only to prove the composition is injective.
Take a nonzero element $x = \tau_0(x') \in \Im^{\!0} \tau_0$ such that
$\rho_2(x) \in \Im^{\!0} \rho_2$.
By the first assertion of Claim \ref{LemmaClaim3},
if $\rho_2(x) \neq 0$, there exists $y \in H^2(A)$
such that the sum of the cup product pairings
$\rho_2(x) \cup \rho_2(y)$ is nonzero.
Since $\rho_2,\tau_2$ are dual to each other, we have
$$ 0 \neq \rho_2(x) \cup \rho_2(y)
      = (\tau_2 \circ \rho_2(x)) \cup y
      = (\tau_2 \circ \rho_2 \circ \tau_0(x')) \cup y
      = 0, $$
which is absurd.
Hence $\rho_2(x) = 0$.
By the second assertion of Claim \ref{LemmaClaim3},
there exists $\tau_0(y') \in \Im^{\!0} \tau_0$
such that the sum of the Lefschetz pairings
$L \cup \tau_0(x') \cup \tau_0(y')$ is nonzero.
Since $\rho_4, \tau_0$ are dual to each other, we have
\begin{align*}
 0 &\neq L \cup \tau_0(x') \cup \tau_0(y')
     = \rho_4(L \cup \tau_0(x')) \cup y'
     = L \cup (\rho_2 \circ \tau_0(x')) \cup y' \\
     &= L \cup \rho_2(x) \cup y'
     = 0,
\end{align*}
which is absurd.
Hence we have Claim \ref{LemmaClaim4}.
\end{proof}

By Claim \ref{LemmaClaim4},
the surjection $\Im \rho_2 \to \Im^{\!1} \rho_2$
has a canonical splitting. Therefore,
we have a canonical decomposition of $\Im \rho_2$ as follows
$$ \Im \rho_2 = \Im^{\!0} \rho_2 \oplus \Im^{\!1} \rho_2. $$
We shall show that this decomposition is an orthogonal decomposition
with respect to the sum of the cup product pairings on $H^2(B)$.
Take
$(\rho_2(a), \rho_2(b)) \in \Im^{\!0} \rho_2 \oplus \Im^{\!1} \rho_2$.
By Claim \ref{LemmaClaim4},
if we choose an appropriate $b$, we can write
$b = \tau_0(c)$ for some $c \in H^0(B)$.
Then, we have
$$ \rho_2(a) \cup \rho_2(b) = \rho_2(a) \cup (\rho_2 \circ \tau_0(c))
       = a \cup (\tau_2 \circ \rho_2 \circ \tau_0(c)) = 0 $$
since $\rho_2,\tau_2$ are dual to each other
and $\tau_2 \circ \rho_2 \circ \tau_0 = 0$.

Now we shall prove Lemma \ref{KeyLemma}.
We have only to prove
$\Ker \tau_2 \cap \Im \rho_2 \subset \Im(\rho_2 \circ \tau_0)$
since the opposite inclusion is trivial.
Let $x \in \Ker \tau_2 \cap \Im \rho_2$.
Then, for all $y \in H^2(A)$,
$$ x \cup \rho_2(y) = \tau_2(x) \cup y = 0 $$
because $\tau_2(x) = 0$.
We write $x = a + b,\ a \in \Im^{\!0} \rho_2,\ b \in \Im^{\!1} \rho_2$.
If $a \neq 0$, by Claim \ref{LemmaClaim3},
there exists $c \in \Im^{\!0} \rho_2$ such that
$$ x \cup c = a \cup c + b \cup c = a \cup c \neq 0, $$
which is absurd.
Hence $a = 0$ and
$ x \in \rho_2(\Im^{\!0} \tau_0) \subset \Im (\rho_2 \circ \tau_0). $
This proves Lemma \ref{KeyLemma}.

Therefore the proof of Proposition \ref{WMC3fold}
and hence Theorem \ref{MainTheorem} is completed.
\end{proof}

\begin{rem}
Lemma \ref{KeyLemma}
is an analogue of Lemma 4.1.9 in \cite{MSaito1}
which is the key step in the proof of
a Hodge analogue of Conjecture \ref{WMC_semistable}.
Note that in the proof of Lemma \ref{KeyLemma},
Claim \ref{LemmaClaim3} is the only part
where we use the conditions of Proposition \ref{WMC3fold}
and the Hodge index theorem.
\end{rem}

\section{A $p$-adic analogue}
\label{SectionPadic}

In this section, we prove a $p$-adic analogue of
Theorem \ref{MainTheorem}
by using the weight spectral sequence of Mokrane
for $p$-adic cohomology.

Let notation be as in \S \ref{SectionIntroduction}.
Assume that $F$ is a finite field.
Let $W(F)$ be the ring of Witt vectors with coefficients in $F$,
and $K_0$ the field of fractions of $W(F)$.
Mokrane constructed a $p$-adic analogue of
the weight spectral sequence of Rapoport-Zink
of the following form
\begin{center}
\begin{tabular}{l}
$\displaystyle E_1^{-r,\,w+r} = \bigoplus_{k \geq \max\{0,-r\}}
      H^{w-r-2k}_{\text{crys}}(X^{(2k+r+1)}/W(F))(-r-k)$ \\
\hspace*{3in} $\displaystyle \Longrightarrow
\quad H^w_{\text{log-crys}}(X_F^{\times}/W(F)^{\times}),$
\end{tabular}
\end{center}
where $H^{\ast}_{\text{crys}}$
denotes the crystalline cohomology,
and $H^w_{\text{log-crys}}(X_F^{\times}/W(F)^{\times})$
denotes the log crystalline cohomology of the special fiber
$X_F = \X \otimes_{\O_K} F$ endowed with a natural log structure
(see \cite{Mokrane}, \S 3.23, Th\'eor\`eme 3.32).

This spectral sequence has similar properties as the $l$-adic case
(compare with Remark \ref{PropertiesWeightSpectralSequence}).
This spectral sequence degenerates at $E_2$
modulo torsion, which is a consequence of the Weil conjecture
for crystalline cohomology
proved by Katz-Messing (\cite{KatzMessing}).
There is a monodromy operator $N$
satisfying the same properties as the $l$-adic case
($N$ coincides with $\nu$ in \cite{Mokrane}, \S 3.33).
Moreover, there is a $p$-adic analogue of
the weight-monodromy conjecture
(\cite{Mokrane}, Conjecture 3.27, \S 3.33).

We have a $p$-adic analogue of Theorem \ref{MainTheorem}.

\begin{thm}
\label{MainTheoremPadic}
Let $X$ be a proper smooth variety of dimension 3 over $K$
which has a proper strictly semistable model $\X$ over $\O_K$.
Let $X_1,\ldots,X_m$ be the irreducible components
of the special fiber of $\X$.
Assume that the following conditions are satisfied.
\begin{enumerate}
\item $\X$ is projective over $\O_K$.
\item $\rho((X_i)_{\overline{F}}) =
  \dim_{K_0} H^2_{\text{\rm crys}}(X_i/W(F)) \otimes_{W(F)} K_0$
  for all $i$.
\end{enumerate}
Then, $N^r$ induces an isomorphism
$$
\begin{CD}
N^r : E_2^{-r,\,w+r}(r) \otimes_{W(F)} K_0
  @>{\cong}>> E_2^{r,\,w-r} \otimes_{W(F)} K_0
\end{CD}
$$
for all $r,w$.
Therefore, a $p$-adic analogue of the weight-monodromy
conjecture (\cite{Mokrane}, Conjecture 3.27, \S 3.33)
holds for $\X$.
\end{thm}

\begin{proof}
The proof is the same as the $l$-adic case.
Namely, for $H^0, H^1$ and the Poincar\'e dual of them,
we can use the same argument as in
Mokrane's proof of a $p$-adic analogue of
the weight-monodromy conjecture for curves and surfaces
(\cite{Mokrane}, \S 5, \S 6).
To prove an analogue of Lemma \ref{KeyLemma},
we use the cycle map of crystalline cohomology
instead of $l$-adic cohomology
(\cite{GilletMessing}, \cite{Gros}).
\end{proof}

\begin{rem}
It is known that
$\dim_{K_0} H^2_{\text{crys}}(X_i/W(F)) \otimes_{W(F)} K_0$
is equal to
$\dim_{\Q_l} H^2((X_i)_{\overline{F}},\Q_l)$
for all $l \neq \text{\rm char}\,F$
(\cite{KatzMessing}, Corollary 1).
Therefore, the conditions of Theorem \ref{MainTheorem} and
Theorem \ref{MainTheoremPadic} are equivalent to each other.
\end{rem}

\end{document}